\documentclass{amsart}

\usepackage{amssymb}
\usepackage{amscd}
\usepackage{latexsym}
\usepackage{amsfonts}
\usepackage{amsmath}
\usepackage{ragged2e}
\usepackage[latin1]{inputenc}

\newtheorem{proposition}{Proposition}
\newtheorem{lemma}{Lemma}
\newtheorem{theorem}{Theorem}
\newtheorem{definition}{Definition}
\newtheorem{remark}{Remark}

\newcommand{\eproof}{\begin{flushright} $\square$ \end{flushright}}

\DeclareMathOperator\sll{sl}

\newcommand{\im}{\mathop{\fam0 Im}\nolimits}

\newcommand{\tr}{\mathop{\fam0 Tr}\nolimits}

\newcommand{\rank}{\mathop{\fam0 Rank}\nolimits}

\newcommand{\Hom}{\mathop{\fam0 Hom}\nolimits}
\newcommand{\End}{\mathop{\fam0 End}\nolimits}

\newcommand{\ad}{\mathop{\fam0 ad}\nolimits}

\newcommand{\Id}{\mathop{\fam0 Id}\nolimits}

\newcommand{\bC}{{\mathbb C}}
\newcommand{\bR}{{\mathbb R}}
\newcommand{\bT}{{\bar T}}
\newcommand{\C}{C}
\newcommand{\A}{{\mathcal A}}

\newcommand{\Z}{{\mathbb Z}}
\newcommand{\bZ}{\Z{}}

\newcommand{\D}{{\mathcal D}}
\newcommand{\ra}{\mathop{\fam0 \rightarrow}\nolimits}

\newcommand{\M}{{\mathcal M}}

\newcommand{\dbar}{ \bar \partial}

\newcommand{\cH}{ {\mathcal H}}
\newcommand{\cM}{ {\mathcal M}}
\newcommand{\T}{ {\mathcal T}}

\renewcommand{\ad}{\mathop{\fam0 ad}\nolimits}

\renewcommand{\L}{{\mathcal L}}

\newcommand{\tD}{ {\tilde D}}
\newcommand{\tG}{ {\tilde G}}
\newcommand{\tW}{ {\tilde W}}
\newcommand{\s}{\sigma}
\newcommand{\bN}{{\mathbb N}}
\newcommand{\trho}{{\tilde \rho}}
\newcommand{\Nabla}{{\mathbf {\hat \nabla}}}
\newcommand{\Nablat}{{\mathbf {\hat \nabla}}^{t}}
\newcommand{\Nablae}{{\mathbf {\hat \nabla}}^e}
\newcommand{\Nablaet}{{\mathbf {\hat \nabla}}^{e,t}}
\newcommand{\BTstar}{\star^{\text{\tiny BT}}}
\newcommand{\tBTstar}{{\tilde \star}^{\text{\tiny BT}}}

\newcommand{\Ric}{\mathop{\fam0 Ric}\nolimits}

\begin{document}

\title[Hitchin's connection and invariant deformation quantization]{Hitchin's connection,
 Toeplitz operators and symmetry invariant deformation quantization}

\author{J{\o}rgen Ellegaard Andersen}
\address{Department of Mathematics\\
        University of Aarhus\\
        DK-8000, Denmark}
\email{andersen@imf.au.dk}

\begin{abstract}

We establish that Hitchin's connection exist for any rigid
holomorphic
family of K{\"a}hler structures on any compact pre-quantizable
symplectic manifold which satisfies certain simple topological
constraints. Using Toeplitz operators we prove that Hitchin's
connection induces a unique formal connection on smooth functions on
the symplectic manifold. Parallel transport of this formal
connection produces equivalences between the corresponding
Berezin-Toeplitz deformation quantizations. - In the cases where the
Hitchin connection is projectively flat, the formal connections will
be flat and we get a symmetry-invariant formal quantization. If a
certain cohomological condition is satisfied a global trivialization
of this algebra bundle is constructed. As a corollary we get a
symmetry-invariant deformation quantization.

Finally, these results are applied to the moduli space situation
in which Hitchin originally constructed his connection. First
we get a proof that the Hitchin connection in this case is the same as the connection
constructed by Axelrod, Della Pietra and Witten. Second we
obtain in this way a {\em mapping class group invariant} formal quantization
of the smooth symplectic leaves of the moduli space
of flat $SU(n)$-connections on any compact surface.
\end{abstract}

\maketitle

\section{Introduction}

In \cite{H} Hitchin introduced a connection over Teichmüller space in the bundle one obtains by applying
geometric quantization to the moduli spaces of flat $SU(n)$-connections. Furthermore Hitchin proved that this connection
is projectively flat. Hitchin's construction was motivated by Witten's study of quantum Chern-Simons theory in $2+1$ dimensions in \cite{W1}.
In fact, Witten constructed via path integral techniques a quantization of
Chern-Simons theory in $2+1$ dimensions and argued in \cite{W1}
that this produced a TQFT indexed by a compact simple Lie group and an
integer level $k$.

Combinatorially this theory was first constructed by Reshetikhin \&
Turaev using representation theory of $U_q(\sll(n,\bC))$ at
$q=e^{(2\pi i)/(k+n)}$ \cite{RT1} and \cite{RT2}. Subsequently these
TQFT's were constructed using skein theory by Blanchet,
Habegger, Masbaum \& Vogel in \cite{BHMV1}, \cite{BHMV2} and
\cite{B1}. In particular these TQFT's provide projective representations of the mapping class groups.
The fact that these representations agree with the representations obtained from the projective
action of the mapping class group on the projectively covariant constant sections of the Hitchin connection
follows by combining the results of a series of papers. First of all, the work of Laszlo \cite{La1} provides an identification
of the Hitchin connection with the TUY-connection constructed in the bundle of WZW-conformal blocks over Teichmüller space
in \cite{TUY}. In joint work with Ueno (\cite{AU1}, \cite{AU2},
\cite{AU3} and \cite{AU4}), we have given a proof, based mainly on
the results of \cite{TUY}, that the TUY-construction of the
WZW-conformal field theory after twist by a fractional power of an
abelian theory, satisfies all the axioms of a modular functor. Furthermore, we have proved that the full $2+1$-dimensional TQFT
that results from this is isomorphic to the one constructed by
BHMV via skein theory as mentioned above.

From Witten's path integral formulation of these theories, one expects that these TQFT's have asymptotic expansions
in the level $k$ of the theory. In \cite{A1} we considered this question in the abelian case and described how one should
make this question precise in the context of the mapping class group representations. As we described in that paper, the
natural asymptotic expansion for these sequences of representations is a mapping class group invariant deformation quantization of the moduli
space of flat $U(1)$-connections. In this paper we extend these results to the case of the above mentioned representations and the $SU(n)$-moduli space (see
Theorem \ref{MCGinvDG} below).
In fact we consider a more general setting, in which we can construct the Hitchin connection, build a formal Hitchin connection and understand
its relation to the associated Berezin-Toeplitz deformation quantizations. Let us describe the generalized setting
we will consider.

Let $(M,\omega)$ be a compact symplectic manifold. Let $I$ be a
family of K\"{a}hler structures on $(M,\omega)$ parameterized
holomorphically
by some complex manifold $\T$. Suppose $V$ is a
vector field on $\T$. Then we can differentiate $I$ along $V$ and we
denote this derivative $V[I] : \T \ra C^\infty(M,\End(TM_\bC))$. We
define $\tG(V) \in C^\infty(M , S^2(TM_\bC))$ by
\[V[I] = \tG(V) \omega,\]
and define $G(V) \in C^\infty(M, S^2(T_\s))$ \footnote{We
denote the holomorphic tangent bundle of $(M,I_\s)$ by $T_\s$.} such
that
\[\tG(V) = G(V) + {\overline G(V)} \]
for all real vector fields $V$ on $\T$. We see that $\tG$ and $G$
are one-forms on $\T$ with values in $C^\infty(M , S^2(TM_\bC))$ and
$C^\infty(M, S^2(T_\s))$ respectively.

\begin{definition}
We say that a complex family $I$ of K\"{a}hler structures on
$(M,\omega)$ is {\em Rigid} if
\[\dbar_\sigma (G(V)_\sigma) = 0 \]
for all vector fields $V$ on $\T$ and all points $\sigma\in \T$.
\end{definition}

Assume now that $(M,\omega)$ is prequantizable. That means there
exists a Hermitian line bundle $(\L, (\cdot,\cdot))$ over $M$ with a
compatible connection $\nabla$ such that
\[F_\nabla = \frac{i}{2\pi}\omega.\]

For every $\sigma\in \T$ we consider the finite dimensional subspace
of $C^\infty(M,\L^k)$ given by
\[H_\sigma^{(k)} = H^0(M_\s, \L^k) = \{s\in C^\infty(M, \L^k)| \nabla^{0,1}_\s s =0 \}.\]
We will assume that these subspaces of holomorphic sections form a
smooth finite rank subbundle $H^{(k)}$ of the trivial bundle
$\cH^{(k)} = \T\times C^\infty(M,\L^k)$.

\begin{theorem}\label{MainGHCI}
Suppose that $I$ is a rigid family of K\"{a}hler
structures on the symplectic prequantizable compact manifold
$(M,\omega)$, which satisfies that there exist $n\in \bZ$ such that
the first Chern class of $(M,\omega)$ is $n [\omega]\in H^2(M,\bZ)$
and $H^1(M,\bR) = 0$. Then the Hitchin connection $\Nabla$ in the
bundle $\cH^{(k)}$ preserves the subbundle $H^{(k)}$. It is given by
\[\Nabla_V = \Nablat_V - u(V)\]
where $\Nablat$ is the trivial connection in $\cH^{(k)}$, $V$ is any smooth vector field
on $\T$ and
$u(V)$ is the second order differential operator given by
\begin{equation}
u(V)(s) = \frac1{2k+n} \left\{ \frac12 \Delta_{G(V)}(s) -  \nabla_{G(V)dF}(s) + 2k
V'[F]s\right\},\label{Hitcon}
\end{equation}
where $V'$ denotes the $(1,0)$-part of the vector field $V$ on $\T$
and $F : \T \ra \C^\infty_0(M)$ is determined by $F_\sigma\in
\C^\infty_0(M)$ being the Ricci-potential for $(M,I_\s)$ for all
$\sigma \in \T$.
\end{theorem}

We prove this theorem in section \ref{ghc}. When we apply this
theorem to the gauge theory example discussed in Section \ref{sec5},
we get as a corollary that Hitchin's connection constructed in \cite{H} is the same
as the one constructed by Axelrod, della Pietra and Witten in
\cite{ADW}.

\begin{remark}
In \cite{AGL}, in joint work with Gammelgaard and Lauridsen, we use half-forms and the metaplectic correction to
prove the existence of a Hitchin connection in the context of
half-form quantization. The assumption that the first Chern class
of $(M,\omega)$ is $n [\omega]\in H^2(M,\bZ)$ is then just
replaced by the vanishing of the second Stiefel-Whitney class of
$M$.
\end{remark}

Returning to the situation at hand, suppose $\Gamma$ is a group
which act on $\T$ and on $M$, such that $I$ is $\Gamma$ equivariant.
Assume further that there is an action of $\Gamma$ on the prequantum
bundle $(\L, (\cdot,\cdot), \nabla)$ covering the $\Gamma$-action on
$M$. It then follows that Hitchin's connection is
$\Gamma$-invariant.

Since $H^{(k)}_\s$ is a finite dimensional subspace of
$C^\infty(M,\L^k)= \cH_\s^{(k)}$ and therefore closed, we have the
orthogonal projection $\pi_\s^{(k)} : \cH_\s^{(k)} \ra H^{(k)}_\s$.
Since $H^{(k)}$ is a smooth subbundle of $\cH^{(k)}$ the projections
$\pi_\s^{(k)}$ form a smooth map $\pi^{(k)}$ from $\T$ to the space
of bounded operators on the $L_2$-completion of $C^\infty(M,\L^k)$.

From these projections we can construct the Toeplitz operators
associated to any smooth function $f\in C^\infty(M)$,
$T^{(k)}_{f,\s} : \cH_\s^{(k)} \ra H^{(k)}_\s$, defined by
\[T_{f,\s}^{(k)}(s) = \pi_\s^{(k)}(fs)\]
for any element $s$ in $\cH_\s^{(k)}$ and any point $\s\in \T$. We
observe that the Toeplitz operators are smooth sections
$T_{f}^{(k)}$ of the bundle $\Hom(\cH^{(k)},H^{(k)})$ and restricts
to smooth sections of $\End(H^{(k)})$.

Let $\D(M)$ be the space of smooth differential operators on $M$
acting on smooth functions on $M$. Let $\C_h$ be the trivial
$C^\infty(M)[[h]]$-bundle over $\T$.

\begin{definition}\label{fc}
A formal connection $D$ is a connection in $\C_h$ over $\T$ of the
form
\[D_V f = V[f] + \tD(V)(f),\]
where $\tD$ is a smooth one-form on $\T$ with values in $\D_h(M) =
\D(M)[[h]]$,  $f$ is any smooth section of $\C_h$, $V$ is any smooth
vector field on $\T$ and $V[f]$ is the derivative of $f$ in the
direction of $V$.
\end{definition}

For a formal connection we get the series of differential operators
$\tD^{(l)}$ given by
\[\tD(V) = \sum_{l=0}^\infty \tD^{(l)}(V) h^l.\]

From Hitchin's connection in $H^{(k)}$ we get an induced connection
$\Nablae$ in the endomorphism bundle $\End(H^{(k)})$. The Toeplitz
operators are not covariant constant sections with respect to
$\Nablae$. They are asymptotically in $k$ in the following very
precise sense.
\begin{theorem}\label{MainFGHCI}
There is a unique formal Hitchin connection $D$ which satisfies that
\begin{equation}
\Nablae_V T^{(k)}_f \sim T^{(k)}_{(D_V f)(1/(2k+n))}\label{Tdf}
\end{equation}
for all smooth section $f$ of $\C_h$ and all smooth vector fields on
$\T$. Moreover
\[\tD = 0 \mod h.\]
Here $\sim$ means the following: For all $L\in \Z_+$ we have that
\[\left\| \Nablae_V T^{(k)}_{f} - \left( T^{(k)}_{V[f]} +
\sum_{l=1}^L T_{\tD^{(l)}_V f}^{(k)} \frac1{(2k+n)^{l}}\right)
\right\| = O(k^{-(L+1)})\] uniformly over compact subsets of $\T$
for all smooth maps $f:\T \ra C^\infty(M)$.
\end{theorem}

This theorem is proved in section \ref{fgHc}. For an
explicit formula for $\tD$ see (\ref{formalcon}).

 Again in the presence of a symmetry group $\Gamma$ as
before, the formal Hitchin connection becomes $\Gamma$-invariant.

By the work of Bordeman, Meinrenken and Schlichenmaier
\cite{BMS}, \cite{Sch}, \cite{Sch1} and \cite{Sch2} applied to the K\"{a}hler manifold
$(M,\omega,I_\sigma)$, we know that for any $f,g\in \C^\infty(M)$,
there is an asymptotic expansion
\[T_{f,\s}^{(k)}T_{g,\s}^{(k)} \sim \sum_{l=0}^\infty
T_{c_\s^{(l)}(f,g),\s}^{(k)} k^{-l},\]
where $c_\s^{(l)}(f,g)\in \C^\infty(M)$ are uniquely
determined. Moreover it gives a deformation quantization
\[f\star_\sigma^{BT}g = \sum_{l=0}^\infty c_\s^{(l)}(f,g)h^l\]
which is known as the Berezin-Toeplitz quantization. By the work
of Karabegov and Schlichenmaier \cite{KS}, it is known to be a differential
deformation quantization.

\begin{proposition}\label{deriv}
For every vector field $V$ on $\T$, the formal operator $D_V$ is a derivation for
$\star_\sigma^{BT}$.
\end{proposition}

Let $\A_h$ be the vector space of sections of $\C_h$ which are
covariant constant with respect to the formal Hitchin connection
$D$ over $\T$. Then by Proposition \ref{deriv}, we see that star
products $\star_\sigma^{BT}$, $\sigma\in \T$, induces an
associative algebra structure on $\A_h$. Moreover, the symmetry
group $\Gamma$ will act by automorphisms of $\A_h$.

\begin{theorem}\label{FQ}
If the formal Hitchin connection has
trivial global holonomy over $\T$, then the algebra $\A_h$ is a
formal quantization of the Poisson algebra of smooth functions on
$(M,\omega)$.
\end{theorem}

\begin{remark}
If the Hitchin connection is projectively flat and $\T$ is simply
connected, then the formal Hitchin connection has
trivial global holonomy over $\T$.
\end{remark}

A formal trivialization of a formal connections is defined as
follows.

\begin{definition}\label{formaltrivi}
A formal trivialization of a formal connection $D$ is a smooth map
$P : \T \ra \D_h(M)$ which modulo $h$ is an isomorphism for all
$\s\in \T$ and such that
\[D_V(P(f)) = 0\]
for all vector fields $V$ on $\T$ and all $f\in C^\infty_h(M)$.
\end{definition}

The existence of a formal trivialization is of course equivalent
to the triviality of the global holonomy of $D$ over $\T$.

In section \ref{PTSI} we construct a global $\Gamma$-equivariant formal trivialization
 under the assumption that
the Hitchin connection is flat and that the first
$\Gamma$-equivariant cohomology of $M$ with coefficients in the
$\Gamma$-module consisting of all differential operators on $M$
vanishes. This further leads to the construction of a $\Gamma$
invariant $*$-product on $M$, simply because a global trivialization
induces a vector space isomorphism between $\A_h$ and
$\C_h^\infty(M)$.

\begin{theorem}\label{general}
Assume that the formal Hitchin connection $D$ is flat and
$$H^1_\Gamma(\T,D(M)) = 0,$$ then there is a $\Gamma$-invariant
trivialization $P$ of $D$ and the $*$-product
\[f\star g = P_\s^{-1}(P_\s(f) \BTstar_\s P_\s(g))\]
is independent of $\s\in \T$ and $\Gamma$-invariant.
\end{theorem}

As an application we apply our results to the moduli space of flat
connections on a surface.

Let $\Sigma$ be a compact two dimensional manifold and $\M$ the
moduli space of flat $SU(n)$-connections on $\Sigma$.

Let $\Gamma$ be the mapping class group of $\Sigma$. The
$\Gamma$ action on $\M$ is Poisson and it preserving all the
symplectic leaves of $\M$. Let $\T$ be Teichm\"{u}ller space of
$\Sigma$. Then for each smooth symplectic leaf $(M,\omega)$ of $\M$
we have a holomorphic map $I$ from $\T$ to the space of complex
structures on $(M,\omega)$ which is $\Gamma$-equivariant.

Applying the above to this moduli space situation, we get that

\begin{theorem}\label{MCGinvDG}
There is a mapping class group invariant formal quantization
on the smooth symplectic leaves $(M,\omega)$ of the moduli space
of flat $SU(n)$-connections on any compact surface.
\end{theorem}

It seems a very interesting open problem to understand how this formal quantization is related to the
quantization of the moduli spaces presented in \cite{AMR1} and \cite{AMR2}, as we have done in
the abelian case in \cite{A1}.

We remark that the interplay between Toeplitz operators and
Hitchin's connection forms the foundation for the proof of the asymptotic
faithfulness theorem in \cite{A2}, the determination of the Nielsen-Thurston types of
mapping classes via TQFT \cite{A3} and forms our example of a representation of the mapping class group,
which has no fixed vectors, but which has an almost fixed vector, i.e. proving the mapping class groups do not have
Kazhdan's property (T) \cite{A4}.

\section{The Hitchin connection}\label{ghc}

Let $(M,\omega)$ be a compact symplectic manifold.

\begin{definition}\label{prequantumb}
A prequantum line bundle $(\L, (\cdot,\cdot),
\nabla)$ over the symplectic manifold $(M,\omega)$ consist of a
complex line bundle $\L$ with a Hermitian structure
$(\cdot,\cdot)$ and a compatible connection  $\nabla$
whose curvature is
\[F_\nabla = \frac{i}{2\pi}\omega,\]
e.g.
\[\nabla_X\nabla_Y - \nabla_Y\nabla_X -  \nabla_{[X,Y]} = \omega(X,Y)\]
for all vector fields $X,Y$ on $M$.
We say that the symplectic manifold $(M,\omega)$ is prequantizable if
there exist a prequantum line bundle over it.
\end{definition}

Recall that the condition for the existence of a prequantum line
bundle is that $[\omega]\in \im(H^2(M,\bZ) \ra H^2(M,\bR))$ and that
the inequivalent choices of prequantum line bundles (if they exist)
are parameterized by $H^1(M,U(1))$. See e.g. \cite{Woodhouse}.

We shall assume
that $(M,\omega)$ is prequantizable and fix a prequantum line bundle $(\L, (\cdot,\cdot),
\nabla)$ over $(M,\omega)$.

Assume that $\T$ is a smooth manifold which smoothly parametrizes
K\"{a}hler structures on $(M,\omega)$. This means we have a
smooth\footnote{Here a smooth map from
$\T$ to $C^\infty(M,W)$ for any smooth vector bundle $W$ over $M$
means a smooth section of $\pi_M^*(W)$ over $\T\times M$, where
$\pi_M$ is the projection onto $M$. Likewise a smooth $p$-form on
$\T$ with values in $C^\infty(M,W)$ is by definition a smooth
section of $\pi_{\T}^*\Lambda^p(\T)\otimes \pi_M^*(W)$ over $\T\times
M$. We will also encounter the situation where we have a bundle
$\tW$ over $\T\times M$ and then we will talk about a smooth $p$-form on $\T$ with values in
$C^\infty(M,\tW_\s)$ and mean a smooth section of
 $\pi_{\T}^*\Lambda^p(\T)\otimes \tW$ over $\T\times M$.}
map $I : \T \ra C^\infty(M,\End(TM))$ such that $(M,\omega, I_\s)$
is a K\"{a}hler manifold for each $\s\in \T$.

We will use the notation $M_\sigma$
for the complex manifold $(M, I_\s)$. For each $\s\in \T$ we
use $I_\s$ to split the the complexified tangent bundle $TM_\bC$ into the
holomorphic and the anti-holomorphic parts, which we denote
$$T_{\s} = E(I_\s,i) = \im(\Id - iI_\s)$$
and
$$\bT_{\s}= E(I_\s,-i) = \im(\Id + iI_\s)$$
respectively.

The real K\"{a}hler-metric $g_\s$ on $(M_\s,\omega)$ extended complex linearly to $TM_\bC$
is by definition
\[g_\s(X,Y) = \omega(X,I_\s Y),\]
where $X,Y \in C^\infty(M,TM\otimes \bC)$.
Both $g_\s$ and
$\omega$ induce isomorphisms $i_{g_\s}, i_\omega : TM_\bC \ra
T^*M_\bC$ and they are related by
\[i_{g_\s} = - I_\s i_\omega. \]
We record for later use that since $\Lambda^2 (I_\s) \omega = \omega$,
we get that
\begin{equation}
(\Lambda^2 i_{g_\s})^{-1}(\omega) = (\Lambda^2 i_{\omega})^{-1}(\omega).\label{metiso=sympiso}
\end{equation}

Suppose $V$ is a vector field on $\T$. Then we can differentiate
$I$ along $V$ and we denote this derivative $V[I] : \T \ra
C^\infty(M,\End(TM_\bC))$. Differentiating the equation $I^2 = -\Id$,
we see that $V[I]$ anti-commutes with $I$. Hence we get that
\[V[I]_\s \in C^\infty(M, (T_\s^*\otimes \bT_\s)\oplus (\bT_\s^*\otimes T_\s))\]
for each $\s\in \T$. Let
$$V[I]_\s = V[I]''_\s + V[I]'_\s $$
be the corresponding decomposition such that $V[I]''_\s\in
C^\infty(M, T_\s^*\otimes \bT_\s)$ and $V[I]'_\s\in
C^\infty(M, \bT_\s^*\otimes T_\s)$.

Now we will further assume that $\T$ is a complex manifold and that $I$ is
a holomorphic
map from $\T$ to the space of all complex structures on $M$.
Concretely, this means that
\[V'[I]_\s = V[I]'_\s\]
and
\[V''[I]_\s = V[I]''_\s\]
for all $\s\in \T$, where $V'$ means the $(1,0)$-part of $V$ and $V''$ means the $(0,1)$-part of
$V$ over $\T$.

We see that
\[V[g](X,Y) = \omega (X, V[I] Y).\]
Since $\omega$ is of type $(1,1)$ and $g$ is symmetric, we see
that
\[V[g] \in C^\infty(M, S^2(T^*) \oplus S^2(\bar T^*))\]
and self-conjugate for real vector fields.
Let us now define $\tG(V) \in C^\infty(M , TM_\bC\otimes TM_\bC)$
by
\[V[I] = \tG(V) \omega,\]
and define $G(V) \in C^\infty(M, T_\s \otimes T_\s)$ such that
\[\tG(V) = G(V) + {\overline G(V)} \]
for all real vector fields $V$ on $\T$. We see that $\tG$ and $G$ are
one-forms on $\T$ with values in $C^\infty(M , TM_\bC\otimes
TM_\bC)$ and $C^\infty(M, T_\s \otimes T_\s)$ respectively.
We observe that
\[V'[I] = G(V)\omega,\]
and $G(V) = G(V')$.

Since $V[g] = \omega
V[I]$, we see that
\[V[g] = \omega \tG(V) \omega = i_\omega\otimes i_\omega (\tG(V)).\]
From this it is clear that $\tG$
takes values in $C^\infty(M , S^2(TM_\bC))$ and therefore that $G$
takes values in $C^\infty(M, S^2(T_\s))$.

On $\L^k$ we have the the smooth family of $\bar\partial$-operators
$\nabla^{0,1}$ defined at $\s\in \T$ by
\[\nabla^{0,1}_\s = \frac12 (1+i I_\s)\nabla.\]

For every $\sigma\in \T$ we consider the finite dimensional
subspace of $C^\infty(M,\L^k)$ given by
\[H_\sigma^{(k)} = H^0(M_\s, \L^k) = \{s\in C^\infty(M, \L^k)| \nabla^{0,1}_\s s =0 \}.\]
We will assume that these subspaces of holomorphic sections form a smooth finite rank
subbundle $H^{(k)}$ of the trivial bundle $\cH^{(k)} = \T\times C^\infty(M,\L^k)$.

Let $\Nablat$ denote the trivial connection in the trivial bundle
$\T\times C^\infty(M,\L^k)$. Let $\D(M,\L^k)$ denote the vector
space of differential operators on $C^\infty(M,\L^k)$. For any
smooth one-form $u$ on $\T$ with values in $\D(M,\L^k)$ we have a
connection $\Nabla$ in $\cH^{(k)}$ given by
$$\Nabla_V = \Nablat_V - u(V)$$
for any vector field $V$ on $\T$.

\begin{lemma}
The connection $\Nabla$ in $\cH^{(k)}$ induces a connection in $H^{(k)}$ if and only if
\begin{equation}
 \frac{i}2  V[I] \nabla^{1,0} s + \nabla^{0,1}u(V)s = 0\label{eqcond}
\end{equation}
for all vector fields $V$ on $\T$ and all smooth sections $s$ of
$H^{(k)}$.
\end{lemma}

\proof
Let $s$ be a section of $H^{(k)}$ over $\T$ and $V$ a vector field
on $\T$. Then $\Nabla_V s$ is a section of
$\cH^{(k)}$, and we compute at a point $\s\in \T$ that
\begin{eqnarray*}
\nabla^{0,1}_\s ((\Nabla_V(s))_\s) & = &
 \nabla^{0,1}_\s (V[s]_\s) - \nabla^{0,1}_\s ((u(V)s)_\s)\\
& = & - \frac{i}2  (V[I] \nabla^{1,0} s)_\s - \nabla^{0,1}_\s ((u(V)s)_\s),
\end{eqnarray*}
since
\[\frac{i}2 (V[I] \nabla^{1,0} s)_\s + \nabla^{0,1}_\s(V[s]_\s) = 0.\]
Hence $\Nabla$ preserves $H^{(k)}$ if and only if \eqref{eqcond}
holds.
\eproof

We observe that
\[V''[I] \nabla^{1,0} s = 0,\]
so $u(V'') = 0$ solves \eqref{eqcond} along the anti-holomorphic
directions on $\T$. In other words the $(0,1)$-part of the trivial
connection $\Nablat$ induces a $\bar\partial$-operator on $H^{(k)}$
and hence makes it a holomorphic vector bundle over $\T$.

This is of course not in general the situation in the $(1,0)$
direction. Let us now consider a particular $u$ and prove that it
solves \eqref{eqcond} under certain conditions.

On the K\"{a}hler manifold $(M_\s,\omega)$ we have the
K\"{a}hler metric and we have the Levi-Civita connection $\nabla$
in $T_\s$. We also have the Ricci potential $F_\s\in
C^\infty_0(M,\bR)$. Here
\[C^\infty_0(M,\bR) = \left\{ f\in C^\infty(M,\bR) \mid \int_M f \omega^m = 0\right\} \]
and the Ricci potential is the element of
$C^\infty_0(M,\bR)$ which satisfies
\[\Ric_\s = \Ric_\s^H + 2 i \partial_\s\dbar_\s F_\s,\]
where $\Ric_\s\in \Omega^{1,1}(M_\s)$ is the Ricci form and $\Ric_\s^H$ is its
harmonic
part. We see that we get this way a smooth function $F : \T \ra C^\infty_0(M,\bR)$.

For any $G\in C^\infty(M, S^2(T_\s))$ we get a linear bundle map
\[G : T^*_\s \ra  T_\s\]
and we have the operator
\begin{eqnarray*}
\Delta_G &: &C^\infty(M,\L^k) \xrightarrow{\nabla^{1,0}_\s} C^\infty(M,T^*_\s\otimes\L^k)
\xrightarrow{G\otimes\Id}
 C^\infty(M,T_\s\otimes\L^k) \\
 && \qquad \xrightarrow{\nabla^{1,0}_\s\otimes \Id +
\Id\otimes \nabla^{1,0}_\s}
 C^\infty(M,T^*_\s\otimes T_\s\otimes\L^k)
\stackrel{\tr}{\ra}C^\infty(M,\L^k).
\end{eqnarray*}

For any smooth function $f$ on $M$, we get a vector field
\[G d f \in C^\infty(M,T_\s).\]
Implicit in this definition is the projection from $TM \cong T_\s
\oplus \bar T_\s$ to $T_\s$, which takes $df$ to $\partial_\s f$.

Putting these constructions together we consider the following
operator for some $n\in \bZ$ such that $2k+n \neq 0$
\begin{equation}
u(V) = \frac1{2k+n}o(V) + V'[F]\label{equ}
\end{equation}
where
\begin{equation}
o(V) = \frac12 \Delta_{G(V)} -  \nabla_{G(V)dF} - n V'[F].\label{eqo}
\end{equation}

The connection associated to this $u$ is denoted $\Nabla$ and we call it the
{\em Hitchin connection}.

\begin{definition}
We say that the complex family $I$ of K\"{a}hler structures on
$(M,\omega)$ is {\em Rigid} if
\[\dbar_\sigma (G(V)_\sigma) = 0 \]
for all vector fields $V$ on $\T$ and all points $\sigma\in \T$.
\end{definition}

We will assume our holomorphic family $I$ is rigid.

\begin{theorem}\label{HCE}
Suppose that $I$ is a rigid family of K\"{a}hler structures
on the compact symplectic prequantizable manifold $(M,\omega)$, which satisfies
that there exist $n\in \bZ$ such that the first Chern class of $(M,\omega)$
is $n [\omega]\in H^2(M,\bZ)$ and $H^1(M,\bR) = 0$. Then $u$ given
by \eqref{equ} and \eqref{eqo} satisfies \eqref{eqcond} for all $k$ such that $2k+n \neq
0$.
\end{theorem}

Hence the Hitchin connection $\Nabla$ preserves the subbundle
$H^{(k)}$ under the stated conditions and we have obtained Theorem \ref{MainGHCI}.

Theorem \ref{HCE} is established through the following three Lemma's.

\begin{lemma} \label{dbarl}
Assume that the first Chern class of $(M,\omega)$
is $n [\omega]\in H^2(M,\bZ)$. For any $\s\in \T$ and
for any $G\in H^0(M_\s, S^2(T_\s))$ we have the following formula
\begin{eqnarray*}
\nabla^{0,1}_\s (\Delta_G(s) -2 \nabla_{G d F_\s}(s)) &=& -
i (2 k + n) G \omega \nabla_\s^{1,0} (s)\\
&  & -ik \tr( - 2  G \partial_\s F \omega + \nabla^{1,0}_\s(G)\omega)
s,
\end{eqnarray*}
for all $s\in H^0(M_\s, \L^k)$.
\end{lemma}

\proof

We compute that
\begin{eqnarray*}
\nabla^{0,1}_\s (\Delta_G(s)) & = & \tr (\nabla^{0,1}\nabla^{1,0}G
\nabla^{1,0}(s))\\
& = & \tr (\nabla^{1,0}\nabla^{0,1}G
\nabla^{1,0}(s))\\
&  & -ik \tr(\omega G \nabla^{1,0}(s)) - i \tr(\Ric_\s G
\nabla^{1,0}(s))\\
& = & -ik \tr(\omega G \nabla^{1,0}(s)) - i \tr(\Ric_\s G
\nabla^{1,0}(s))\\
&  &  -i k \tr(\nabla^{1,0}(G \omega s))\\
& = & -2ik \tr(\omega G \nabla^{1,0}(s)) - i \tr(\Ric_\s G
\nabla^{1,0}(s))\\
&  &  -i k \tr(\nabla^{1,0}(G) \omega) \otimes s,\\
\end{eqnarray*}
since $G$ is holomorphic and $\nabla(\omega) = 0$ because
$(M_\s,\omega)$ is K\"{a}hler. The assumption $c_1(M,\omega) = n
[\omega]$ implies that $\Ric_\s^H= n\omega$ and so
\[\Ric_\s = n\omega + 2 i \partial_\s \bar \partial_\s F_\s.\]
From this we conclude the stated formula.
\eproof

\begin{lemma}\label{Vriccipot}
We have the following relation
$$
2 i \bar\partial_\s (V'[F]_\s) =  \frac12 \tr(  2 G(V)\partial (F) \omega -
\nabla^{1,0}(G(V))\omega)_\s
$$
provided $H^1(M,\bR) = 0$.
\end{lemma}

To prove this lemma, we need a formula for the
variation of the Ricci-form.

\begin{lemma}
For any smooth vector field $V$ on $\T$ we have that
\begin{equation}
(V'[\Ric])^{1,1} = - \frac12 \partial \tr(\nabla^{1,0}(G(V)) \omega).
\end{equation}
\end{lemma}

\proof
The Ricci form $\Ric_\s$ of the K\"{a}hler manifold $(M_\s,\omega)$
is by definition
\[\Ric_\s = R_\s(\omega),\]
where $R_\s\in C^\infty(M, S^2(\Lambda^{1,1}_\s M))$ is the K\"{a}hler
curvature and
$$R_\s(\omega) = R_\s((\Lambda^2 i_{g_\s})^{-1}(\omega)).$$
From this we conclude that
\[V'[\Ric] = V'[R]((\Lambda^2 i_{\omega})^{-1}(\omega))\]
where we have used \eqref{metiso=sympiso}. According to Theorem
1.174, (c), in \cite{Besse}, we have for any four vector fields $X,Y,Z,U$ on $M$, that
\begin{eqnarray*}
V'[R](X,Y,Z,U) & = & \frac12 ( \nabla^2_{Y,Z}(V'[g])(X,U) +
\nabla^2_{X,U}(V'[g])(Y,Z)\\
 & & - \nabla^2_{X,Z}(V'[g])(Y,U) -
\nabla^2_{Y,U}(V'[g])(X,Z)\\
&& + V'[g](R(X,Y)Z,U) - V'[g](R(X,Y)U,Z)),
\end{eqnarray*}
where $R(X,Y) \in C^\infty(M,\End(TM))$ is $(3,1)$-curvature of
$\nabla$, the Levi-Civita connection on $(M_\s,\omega)$
and $\nabla_{X,Y}^2 = \nabla_X \nabla_Y- \nabla_{\nabla_X Y}$.
The Levi-Civita connection $\nabla$
of $(M_\s,\omega)$ and therefore also $R(X,Y)_\s$ preserves the
type decomposition of tensors on $(M_\s,\omega)$. Hence when we
want to compute $V'[R]$ applied to $(\Lambda^2 i_{\omega})^{-1}(\omega)
\in C^\infty(M, T_\s\wedge \bar T_\s)$, and we are only interested in the $(1,1)$
part of the result, we only get contributions from the third and fourth
term. These two terms give
\begin{eqnarray*}(V'[\Ric])^{1,1} & = & - \frac12 ( \tr(\nabla^{1,0}\nabla^{1,0}(V'[g])\otimes
(\Lambda^2 i_{\omega})^{-1}(\omega) ) \\
& & +  \tr(\nabla^{0,1}\nabla^{0,1}(V'[g])\otimes
(\Lambda^2 i_{\omega})^{-1}(\omega))).
\end{eqnarray*}
Using that $\nabla(\omega) = 0 $ and that $V'[g] = i_\omega \otimes i_\omega
G(V)$ we get that
\[(V'[\Ric])^{1,1} = - \frac12 \partial \tr(\nabla^{1,0}(G(V)) \otimes \omega). \]
\eproof

\paragraph{{\bf Proof of Lemma \ref{Vriccipot}.}}
By the definition of the Ricci potential
\[\Ric = \Ric^H + 2 i d \bar \partial F\]
where $\Ric^H = n\omega$ by the assumption. Hence
\[V'[\Ric] = - d V'[I] d F + 2i d \bar \partial V'[F]\]
and therefore
\[2 i \partial \bar \partial V'[F] = (V'[\Ric])^{1,1}
+ \partial V'[I] \partial F.\]
From the above we conclude that
$$\frac12 \tr(2 G(V) \partial F \omega -  \nabla^{1,0}(G(V)) \omega)_\s -
2i \bar \partial_\s V'[F]_\s \in \Omega^{0,1}_\s(M)$$
is $\partial_\s$-closed. By Lemma \ref{dbarl} it is also $\bar\partial_\s$-closed, hence it is a closed one form
on $M$. But since we
assume that $H^1(M,\bR) = 0$, we see it is exact, but then it in fact vanishes since
it is also of type $(0,1)$ on $M_\s$.

\eproof

From the above we conclude that
$$
u(V) = \frac1{2k+n} \left\{ \frac12 \Delta_{G(V)} - \nabla_{G(V)dF} + 2k V'[F]\right\}
$$
solves \eqref{eqcond} and hence we have established
Theorem \ref{HCE}.

We are actually interested in the induced
connection $\Nablae$ in the endomorphism bundle $\End(H^{(k)})$.
Suppose $\Phi$ is a section of $\End(H^{(k)})$. Then for all sections
$s$ of $H^{(k)}$ and all vector fields $V$ on $\T$, we have that
\[(\Nablae_V \Phi) (s) = \Nabla_V \Phi(s) - \Phi(\Nabla_V(s)).\]
Assume now that we have extended $\Phi$ to a section of $\Hom (\cH^{(k)},H^{(k)})$ over $\T$.
Then
\begin{equation}\label{endocon}
\Nablae_V \Phi = \Nablaet_V \Phi + [\Phi, u(V)]
\end{equation}
where $\Nablaet$ is the trivial connection in the trivial
bundle $\End(\cH^{(k)})$ over $\T$.

\begin{lemma}\label{localform}
There exist smooth one-forms $X_r,  Z$ and functions $Y_r,$, $r = 1, \ldots, R$,
on $\T$ with values in $C^\infty(M,T_\s)$  such that
\begin{equation}
\frac12 \Delta_{G(V)} - \nabla_{G(V)d F} = \sum_{r=1}^R \nabla_{X_r(V)}
\nabla_{Y_r} + \nabla_{Z(V)}.\label{Opform}
\end{equation}
for all vector fields $V$ on $\T$.
\end{lemma}

\proof
We fix a finite partition of unity $(U_i, \rho_i)$ of $M$, such that each
$U_i$ is a contractible coordinate neighbourhood $(U_i, x_i)$. Let $\trho_i$
be a smooth function with compact support in $U_i$ which is constant 1 on
the support of $\rho_i$. Using
the coordinate $x_i$  we get a trivialization of $T(U_i)$. By
combining these with the projections onto the varying holomorphic
tangent sub-bundle, we get a smoothly varying trivialization of
these. Using this we see that we for each $i$ can find smooth one
forms $X_i^{(j)}$ and $Y_i^{(j)}$ on $\T$ with values in
$C^\infty(U_i, T(U_i)\otimes \bC)$, which maps $T_\s\T$
to $C^\infty(M,T_\s)$ for all $\s\in\T$ and such that
$$G(V)|_{U_i} = 2 \sum_j X_i^{(j)}(V) Y_i^{(j)}$$
and hence
$$G(V) = 2 \sum_i \sum_j \rho_i X_i^{(j)}(V) \trho_i Y_i^{(j)}.$$
Thus we conclude there exist smooth one-forms $X_r$ and functions
$Y_r$, $r = 1, \ldots, R$, on $\T$ with values in $C^\infty(M,
T_\s)$ such that
\[G(V) = 2 \sum_{r=1}^R X_r(V) Y_r.\]
We now compute that
\begin{eqnarray*}
\frac12 \Delta_{G(V)} -  \nabla_{G(V)dF} & = & \sum_{r=1}^R \nabla_{X_r(V)}
\nabla_{Y_r}\\
& & + \sum_{r=1}^R \tr(\nabla(X_r(V)))
\nabla_{Y_r} - 2 \sum_{r=1}^R X_r(V)(F)
\nabla_{Y_r}.
\end{eqnarray*}
From this we see that
\[Z(V) = \sum_{r=1}^R \tr(\nabla(X_r(V)))
Y_r - 2 \sum_{r=1}^R X_r(V)(F)
Y_r.\]
\eproof

This gives us the expression
\begin{equation}
u(V) = \frac1{2k+n}\left( \sum_{r=1}^R \nabla_{X_r(V)}
\nabla_{Y_r} + \nabla_{Z(V)} - n V'[F]\right)  + V'[F].
\label{Hitchcons}
\end{equation}
All we need to use about $F : \T \ra C^\infty(M)$ below is that it is a
smooth function, such that $F_\sigma$ is real valued
on $M$ for all $\sigma\in \T$.

Suppose $\Gamma$ is a group which acts by bundle automorphisms of
$\L$ over $M$ preserving both the Hermitian structure
and the connection in $\L$. Then there is an induced action of
$\Gamma$ on $(M,\omega)$. We will further assume that $\Gamma$
acts on $\T$ and that $I$ is $\Gamma$-equivariant. In this case we
immediately get the following invariance.

\begin{lemma}
The natural induced action of $\Gamma$ on $\cH^{(k)}$ preserves
the subbundle $H^{(k)}$ and the  Hitchin connection.
\end{lemma}

\section{Berezin-Toeplitz deformation quantization on compact K\"{a}hler manifolds}\label{BZdq}

 For each $f\in C^\infty(M)$ we consider the prequantum operator, namely the
differential operator $P_f^{(k)}: C^\infty(M,L^k) \ra
C^\infty(M,L^k)$ given by
\[P_f^{(k)} = - \frac1k \nabla_{X_f} + i f\cdot\]
where $X_f$ is the Hamiltonian vector field associated to $f$.

These operators act on $C^\infty(M,\L^k)$ and therefore also on the
bundle $\cH^{(k)}$, however, they do not preserve the subbundle
$H^{(k)}$. In order to turn these operators into operators which
act on $H^{(k)}$ we need to consider the Hilbert space structure.

Integrating the inner product of two sections against the volume
form associated to the symplectic form gives the pre-Hilbert space
structure on $C^\infty(M,\L^k)$
\[\langle s_1,s_2\rangle = \frac1{m!}\int_M  (s_1,s_2) \omega^m.\]
We think of this as a pre-Hilbert space structure on the trivial
bundle $\cH^{(k)}$ which of course is compatible with the trivial
connection in this bundle.
This pre-Hilbert space structure induces a Hermitian structure $\langle\cdot,\cdot\rangle$ on
the finite rank subbundle $H^{(k)}$ of $\cH^{(k)}$.
The Hermitian structure $\langle\cdot,\cdot\rangle$ on $H^{(k)}$
also induces the operator norm $\|\cdot\|$ on $\End(H^{(k)})$.

Since $H^{(k)}_\s$ is a finite dimensional subspace
of $C^\infty(M,\L^k)= \cH_\s^{(k)}$ and therefore closed, we have the orthogonal projection
$\pi_\s^{(k)} : \cH_\s^{(k)} \ra
H^{(k)}_\s$. Since $H^{(k)}$ is a smooth subbundle of $\cH^{(k)}$
the projections $\pi_\s^{(k)}$ form a smooth map $\pi^{(k)}$ from $\T$ to the
space of bounded operators on the $L_2$-completion of
$C^\infty(M,\L^k)$. The easiest way to see this is to consider a
local frame for $(s_1, \ldots s_{\rank H^{(k)}})$ of $H^{(k)}$. Let $h_{ij} = \langle s_i,
s_j\rangle$. Let $h^{-1}_{ij}$ be the
inverse matrix of $h_{ij}$. Then
\begin{equation}\label{projf}
\pi_\s^{(k)}(s) = \sum_{i,j}\langle s, (s_i)_\s\rangle (h^{-1}_{ij})_\s (s_j)_\s.
\end{equation}
This formula will be useful when we have to compute the
derivative of $\pi^{(k)}$ along vector fields on $\T$.

From these projections we can construct the Toeplitz operators
associated to any smooth function $f\in C^\infty(M)$, $T^{(k)}_{f,\s} : \cH_\s^{(k)} \ra
H^{(k)}_\s$, defined by
\[T_{f,\s}^{(k)}(s) = \pi_\s^{(k)}(fs)\]
for any element $s$ in $\cH_\s^{(k)}$ and any point $\s\in \T$.
We observe that the Toeplitz operators are smooth sections $T_{f}^{(k)}$ of
the bundle $\Hom(\cH^{(k)},H^{(k)})$ and restrict to smooth sections of $\End(H^{(k)})$.

\begin{remark}
Similarly for any Pseudo-differential operator $A$ on $M$ with
coefficients in $\L^k$ (which may even depend on $\s\in \T$),
we can consider the associated Toeplitz
operator $\pi^{(k)} A$ and think of it as a section of
$\Hom(\cH^{(k)},H^{(k)})$. However, whenever we consider asymptotic expansions of such or
consider their operator norms, we implicitly restrict them to $H^{(k)}$ and consider them as
section of $\End(H^{(k)})$ or equivalently assume that they have been
precomposed with $\pi^{(k)}$.

\end{remark}

We recall by Tuynman's theorem
(see \cite{Tuyn}) that if we compose the prequantum operator
associated to $f$ by the orthogonal projection, then it can be
rewritten as a Toeplitz operator:

\begin{theorem}[Tuynman]\label{Tuynt}
For any $f\in C^\infty(M)$ and any point $\s\in\T$ we have that
\[\pi_\s^{(k)} \circ P_f^{(k)} = i T^{(k)}_{f -\frac1{2k}\Delta_\s f,\s}\]
as operators from $\cH_\s^{(k)}$ to $H^{(k)}_\s$, where $\Delta_\s$ is the Laplacian
on $(M_\s,\omega)$.
\end{theorem}

Tuynman's formula is of course equivalent to
\[\pi^{(k)}_\s\nabla_{X_f} = T^{(k)}_{\frac{i}2 \Delta_\s(f),\s}.\]

This formula is really a corollary of a more
general formula which we will need.

Suppose we have a smooth section $X\in C^\infty(M, T_\s)$ of the holomorphic tangent
bundle of $M_\s$. We then claim that the operator $\pi^{(k)} \nabla_X$
is a zero-order Toeplitz operator. Suppose $s_1\in C^\infty(M,\L^k)$ and $s_2 \in
H^0(M_\s,\L^k)$, then we have that
\[X(s_1,s_2) = (\nabla_X s_1, s_2).\]
Now, calculating the Lie derivative along $X$ of $(s_1,s_2)\omega^m$ and using
 the above,
one obtains after integration that
\[\langle \nabla_X s_1, s_2 \rangle = - \langle \Lambda d(i_X\omega) s_1, s_2 \rangle,\]
where $\Lambda$ denotes contraction with $\omega$. Thus
\begin{equation}\pi^{(k)} \nabla_X = T_{f_X}^{(k)},\label{1to0order}
\end{equation}
as operators from $C^\infty(N,L^k)$ to $H^0(N,L^k)$, where $f_X = -\Lambda d(i_X\omega)$.

Tuynman's formula above now follows from
\[\Lambda d(i_{(X_f)^{(1,0)}}\omega) = -\Lambda d(\bar\partial_\s f) = -
\frac{i}2 \Delta_\s f.\]

Iterating \eqref{1to0order}, we find for all $X_1,X_2 \in C^\infty(M,T_\s)$ that
\begin{equation}\pi^{(k)} \nabla_{X_1}\nabla_{X_2} = T^{(k)}_{f_{X_2}f_{X_1}
- X_2(f_{X_1})}\label{2to0order}
\end{equation}
again as operators from $C^\infty(M,\L^k)$ to $H^0(M_\s,\L^k)$.

For $X\in C^\infty(M, T_\s)$, the complex conjugate vector field
$\bar X \in \C^\infty (M, \bar T_\s)$ is a section of the antiholomorphic tangent bundle,
and for $s_1, s_2 \in C^\infty(M,\L^k)$, we have that
\[\bar X (s_1,s_2) = (\nabla_{\bar X} s_1, s_2) + ( s_1, \nabla_X s_2).\]
Computing the Lie derivative along $\bar X$ of
$(s_1,s_2)\omega^m$ and integrating, we get that
\[\langle \nabla_{\bar X} s_1, s_2 \rangle + \langle (\nabla_{X})^* s_1, s_2 \rangle
 = - \langle \Lambda d(i_{\bar X}\omega) s_1, s_2 \rangle.\]
Hence we see that
\[(\nabla_{X})^* = - \left( \nabla_{\bar X} - f_{\bar X} \right)\]
as operators on $C^\infty(M,\L^k)$. In particular, we see that
\begin{equation}\label{1to0order*}
\pi^{(k)} (\nabla_{X})^* \pi^{(k)} = - T^{(k)}_{f_{\bar X}}|_{H^0(M_\s,\L^k)}
: H^0(M_\s,\L^k) \ra H^0(M_\s,\L^k).
\end{equation}
For two smooth sections $X_1, X_2$ of the holomorphic tangent bundle
$T_\s$ and a smooth function $h\in C^\infty(M)$, we deduce from the
formula for $(\nabla_{X})^*$ that
\begin{eqnarray}\label{2to0order*}
\pi^{(k)} (\nabla_{{X}_1})^*(\nabla_{{X}_2})^* h\pi^{(k)} & = &
\pi^{(k)} \bar X_1 \bar X_2(h)\pi^{(k)} -\\
& & \quad
\pi^{(k)} f_{{\bar X}_1} \bar X_2(h)\pi -\pi f_{{\bar X}_2}\bar X_1(h) \pi^{(k)} -
 \nonumber \\
& & \quad \pi^{(k)}\bar X_1(f_{{\bar X}_2}) h \pi + \pi f_{{\bar X}_1} f_{{\bar
X}_2}h
\pi^{(k)}\nonumber
\end{eqnarray}
as operators on $H^0(M_\s,\L^k)$.

The product of two Toeplitz operators associated to two smooth functions
will in general not be the Toeplitz operator associated to a smooth function again,
but there is an asymptotic expansion of the product in terms of such Toeplitz operators
on a compact K{\"a}hler manifold by the results of Schlichenmaier \cite{Sch}.

\begin{theorem}[Schlichenmaier]\label{S}
For any pair of smooth functions $f_1, f_2\in \C^\infty(M)$, we
have an asymptotic expansion
\[T_{f_1,\s}^{(k)}T_{f_2,\s}^{(k)} \sim \sum_{l=0}^\infty
T_{c_\s^{(l)}(f_1,f_2),\s}^{(k)} k^{-l},\]
where $c_\s^{(l)}(f_1,f_2) \in C^\infty(M)$ are uniquely determined since
$\sim$ means the following: For all $L\in \Z_+$ we have that
\[\|T_{f_1,\s}^{(k)}T_{f_2,\s}^{(k)} - \sum_{l=0}^L T_{c_\s^{(l)}(f_1,f_2),\s}^{(k)} k^{-l}\| =
O(k^{-(L+1)})\]
uniformly over compact subsets of $\T$.
Moreover, $c_\s^{(0)}(f_1,f_2) = f_1f_2$.
\end{theorem}

\begin{remark} It will be useful for us to define new
coefficients ${\tilde c}_\s^{(l)}(f,g) \in C^\infty(M)$ which
correspond to the expansion of the product in $1/(2k+n)$ (where $n$ is some fixed integer):
\[T_{f_1,\s}^{(k)}T_{f_2,\s}^{(k)} \sim \sum_{l=0}^\infty
T_{{\tilde c}_\s^{(l)}(f_1,f_2),\s}^{(k)} (2k+n)^{-l}.\]
\end{remark}

Theorem \ref{S}
is proved in \cite{Sch}, where it is also proved that the formal
generating series for the $c_\s^{(l)}(f_1,f_2)$'s gives a formal
deformation quantization\footnote{We have the opposite sign-convention on the curvature,
which means our $c_l$ are $(-1)^l c_l$ in \cite{Sch}.} of symplectic manifold $(M,\omega)$.

We recall the definition of a formal deformation quantization. Introduce the space of formal
functions $C^\infty_h(M) =
C^\infty(M)[[h]]$ as the space of formal power series in the
variable $h$ with coefficients in $C^\infty(M)$. Let $\bC_h =
\bC[[h]]$.

\begin{definition}
A deformation quantization of $(M,\omega)$ is an associative
product $*$ on $C^\infty_h(M)$ which respects the $\bC_h$-module structure.
It is determined by a sequence of bilinear operators
$$c^{(l)} : C^\infty(M) \otimes C^\infty(M) \ra C^\infty(M)$$
defined through
\[f * g = \sum_{l=0}^\infty c^{(l)}(f,g) h^{l},\]
where $f,g \in C^\infty(M)$. The deformation quantization is said
to be differential, if the operators $c^{(l)}$ are
bidifferential operators. Considering the symplectic action of $\Gamma$ on
$(M,\omega)$, we say that a $*$-product $*$ is $\Gamma$-invariant if
$$ \gamma^*(f*g) = \gamma^*(f)*\gamma^*(g)$$
for all $f,g\in C^\infty(M)$ and all $\gamma\in\Gamma$.
\end{definition}

\begin{theorem}[Karabegov \& Schlichenmaier]\label{tKS1}
The product $\BTstar_\s$ given by
\[f \BTstar_\s g = \sum_{l=0}^\infty (-1)^l c_\s^{(l)}(f,g) h^{l},\]
where $f,g \in C^\infty(M)$ and $c_\s^{(l)}(f,g)$ are
determined by Theorem \ref{S}, is a differentiable deformation quantization of $(M,\omega)$.
\end{theorem}

\begin{definition}
The Berezin-Toeplitz deformation quantization of the compact K{\"a}hler
manifold $(M_\s,\omega)$ is the product $\BTstar_\s$.
\end{definition}

\begin{remark} Let $\Gamma_\s$ be the $\s$-stabilizer subgroup of $\Gamma$.
For any element $\gamma\in \Gamma_\s$, we have that
\[\gamma^*(T^{(k)}_{f,\s}) = T^{(k)}_{\gamma^*f,\s}.\]
This implies the invariance of $\BTstar_\s$ under the $\s$-stabilizer
$\Gamma_\s$.
\end{remark}

\begin{remark} We define a new $*$-product by
\[f\tBTstar_\s g = \sum_{l=0}^\infty (-1)^l {\tilde c}_\s^{(l)}(f,g) h^{l}.\]
Then
\[f\tBTstar_\s g = \left( (f\circ \phi^{-1})\BTstar_\s (g\circ \phi^{-1})\right) \circ \phi\]
for all $f,g\in C_h^\infty(M)$, where $\phi(h) = \frac{h}{2 + nh}$.
\end{remark}

In \cite{KS}, this Berezin-Toeplitz deformation quantization is identified in terms of
Karabegov's classification of $*$-products with separation of variables
on K{\"a}hler manifolds. Adopting the convention where the roles of
holomorphic and anti-holomorphic are interchanged in the condition
for a star product to be with separation of variables
from \cite{KS}, the main result of that paper reads

\begin{theorem}[Karabegov \& Schlichenmaier] \label{tKS2}
The Karabegov form $\tilde \omega_\s$ of the Berezin-Toeplitz $*$-product $\BTstar_\s$ is
\[\tilde \omega_\s = \frac1h \omega + \Ric_\s.\]
\end{theorem}

We will also need the following theorem due to
 Bordemann, Meinrenken and Schlichenmaier (see \cite{BMS}).

\begin{theorem}[Bordemann, Meinrenken and Schlichenmaier]\label{BMS1}
For any $f\in \C^\infty(M)$ we have that
\[\lim_{k\ra \infty}\|T_{f,\s}^{(k)}\| = \sup_{x\in M}|f(x)|.\]
\end{theorem}

Since the association of the sequence of Toeplitz operators
$T^{(k)}_{f,\s}$, $k\in \Z_+$ is linear in $f$, we see from this theorem,
that this association is faithful.

\section{The formal Hitchin connection}\label{fgHc}

We assume the conditions on $(M,\omega)$ and $I$ of Theorem
\ref{HCE}, thus providing us with a  Hitchin connection
$\Nabla$ in $H^{(k)}$ over $\T$ and the associated connection
$\Nablae$ in $\End(H^{(k)})$. Let $\D(M)$ be the space of smooth differential operators on $M$
acting on smooth functions on $M$. Let $\C_h$ be the trivial
$C^\infty_h(M)$-bundle over $\T$.

\begin{definition}\label{fc2}
A formal connection $D$ is a connection in $\C_h$ over $\T$ of the
form
\[D_V f = V[f] + \tD(V)(f),\]
where $\tD$ is a smooth one-form on $\T$ with values in $\D_h(M) =
\D(M)[[h]]$,  $f$ is any smooth section of $\C_h$, $V$ is any smooth
vector field on $\T$ and $V[f]$ is the derivative of $f$ in the
direction of $V$.
\end{definition}

For a formal connection we get the series of differential operators
$\tD^{(l)}$ given by
\[\tD(V) = \sum_{l=0}^\infty \tD^{(l)}(V) h^l.\]

From Hitchin's connection in $H^{(k)}$ we get an induced connection
$\Nablae$ in the endomorphism bundle $\End(H^{(k)})$. The Toeplitz
operators are not covariant constant sections with respect to
$\Nablae$. They are asymptotically in $k$ in the following very
precise sense.
\begin{theorem}\label{MainFGHCI2}
There is a unique formal Hitchin connection $D$ which satisfies that
\begin{equation}
\Nablae_V T^{(k)}_f \sim T^{(k)}_{(D_V f)(1/(2k+n))}\label{Tdf2}
\end{equation}
for all smooth section $f$ of $\C_h$ and all smooth vector fields on
$\T$. Moreover
\[\tD = 0 \mod h.\]
Here $\sim$ means the following: For all $L\in \Z_+$ we have that
\[\left\| \Nablae_V T^{(k)}_{f} - \left( T^{(k)}_{V[f]} +
\sum_{l=1}^L T_{\tD^{(l)}_V f}^{(k)} \frac1{(2k+n)^{l}}\right)
\right\| = O(k^{-(L+1)})\] uniformly over compact subsets of $\T$
for all smooth maps $f:\T \ra C^\infty(M)$.
\end{theorem}

We call this formal connection $D$ the {\em formal Hitchin
connection}.

Before proving this theorem we need to establish some basic
properties.

First, we need a useful formula for the derivative of the
orthogonal projection $\pi^{(k)}$ along a curve $\sigma_t$ in $\T$.
To this end, consider a basis of
covariant constant sections $(s_i)_t$, $i=1,\ldots, \rank H^{(k)}$, of
$H^{(k)}$ over a curve $\sigma_t$ in $\T$:
\[(s_i)'_t = u(\s'_t)((s_i)_t), \qquad i=1, \ldots , \rank H^{(k)}.\]
Recall formula (\ref{projf}) for the projection $\pi^{(k)} : \cH^{(k)}_\s \ra
H^{(k)}_\s$ and compute the derivative along $\s_t$: For
any fixed $s\in C^\infty(M,L^k)$, we have that
\begin{eqnarray*}
(\pi^{(k)}_{\s_t})'(s) & = & \sum_{i,j}\langle s, (s_i)'_{t}\rangle
(h^{-1}_{ij})_{t}(s_j)_{t}\\
&& + \sum_{i,j}\langle s, (s_i)_{t}\rangle (h^{-1}_{ij})'_{t}(s_j)_{t}\\
&& + \sum_{i,j}\langle s, (s_i)_{t}\rangle (h^{-1}_{ij})_{t}(s_j)'_{t}.
\end{eqnarray*}
An easy computation gives that
\[(h^{-1}_{ij})'_{t} = -\sum_{l,r}(h^{-1}_{il})_{t}
(\langle (s_l)'_{t}, (s_r)_{t}\rangle +
\langle (s_l)_{t}, (s_r)'_{t}\rangle) (h^{-1}_{rj})_{t},\]
so
\begin{eqnarray*}
\pi^{(k)}_{\s_t}(\pi^{(k)}_{\s_t})'(s) & = & \sum_{i,j}\langle u_{G(\s'_t)}^* s,
(s_i)_{t}\rangle (h^{-1}_{ij})_{t}(s_j)_{t} \\
&& - \sum_{i,l,m,j}\langle s, (s_i)_{t}\rangle
(h^{-1}_{il})_{t}\langle (s_l)_{t}, (s_m)'_{t}\rangle
(h^{-1}_{mj})_{t}(s_j)_{t}\\
& = & \pi^{(k)}_{\s_t} u(\s'_t)^*(s) - \sum_{m,j}\langle \pi^{(k)}_{\s_t}
s,u(\s'_t)((s_i)_t)\rangle
(h^{-1}_{mj})_{t}(s_j)_{t}\\
& = & \pi^{(k)}_{\s_t}  u(\s'_t)^*(s) - \pi^{(k)}_{\s_t}  u(\s'_t)^* \pi^{(k)}_{\s_t}(s).
\end{eqnarray*}
Hence we conclude that
\begin{lemma} For any smooth vector field $V$ on $\T$, we have
that
\begin{equation}
\pi^{(k)} V[\pi^{(k)}] = \pi^{(k)} u(V)^* -
\pi^{(k)} u(V)^*\pi^{(k)}. \label{derivepi}
\end{equation}
\end{lemma}

Let us now apply this formula for the derivative.

\paragraph{{\bf Proof of Theorem \ref{MainFGHCI}.} }

Using formula \eqref{endocon} we see that
\[\Nablae_V T_f^{(k)} = V[T_f^{(k)}] + [T_f^{(k)}, u(V)]\]
and so we compute
\begin{eqnarray*}
\pi^{(k)} \Nablae_V T_f^{(k)}\pi^{(k)} &=& \pi^{(k)} V[f]\pi^{(k)}
+ \pi^{(k)} V[\pi^{(k)}] f \pi^{(k)}\\
& & - \pi^{(k)} u(V) \pi^{(k)} f \pi^{(k)} + \pi^{(k)} f u(V)
\pi^{(k)}\\
& = & \pi^{(k)} V[f]\pi^{(k)}\\
& &  \pi^{(k)} V''[F] f \pi^{(k)} - \pi^{(k)} V''[F] \pi^{(k)} f
\pi^{(k)}\\
& & + \pi^{(k)} V'[F] f \pi^{(k)} - \pi^{(k)} V'[F] \pi^{(k)} f
\pi^{(k)}\\
& & + \frac1{2k+n}( \pi^{(k)} o(V)^* f \pi^{(k)} - \pi^{(k)} o(V)^* \pi^{(k)} f
\pi^{(k)}\\
& & + \pi^{(k)}  f o(V) \pi^{(k)} - \pi^{(k)} o(V) \pi^{(k)} f
\pi^{(k)})
\end{eqnarray*}
Now by combining \eqref{Hitchcons} with \eqref{1to0order} to \eqref{2to0order*}, we get a one-form $E$ on $\T$ with values in $\D(M)$ such
that
\[T^{(k)}_{E(V)f} = \pi^{(k)} o(V)^* f \pi^{(k)} + \pi^{(k)}  f o(V) \pi^{(k)}. \]
Let $H$ be the one-form on $\T$ with values in $C^\infty(M)$ such
that $H(V) = E(V)(1)$. Then we get the formula
\begin{eqnarray}
\pi^{(k)} \Nablae_V T_f^{(k)}\pi^{(k)} & = & T_{V[f]}^{(k)} +
( T_{V[F]f}^{(k)} - T_{V[F]}^{(k)}T_{f}^{(k)})\nonumber\\
 & & + \frac1{2k+n}(T^{(k)}_{E(V)(f)} -
 T_{H(V)}^{(k)}T_{f}^{(k)}).\label{NablaeT}
 \end{eqnarray}
From this we obtain the wanted estimate by letting
\begin{equation}
\tD(V)(f) =
V[F]f - V[F]\tBTstar f + h (E(V)(f) - H(V) \tBTstar
f)\label{formalcon}
\end{equation}
which is clearly divisible by $h$. \eproof

\begin{lemma} \label{sd=ds}
If $A$ is a smooth family of Toeplitz operators of order $d$, then
$\Nablae_V A$ is also a smooth family of Toeplitz operators of
order $d$ and
\[\s_d(\Nablae_VA) = V[\s_d(A)]\]
for any vector field $V$ on $\T$.
\end{lemma}

\proof
Since $\Nablae_V A$ is again a smooth section of $\End(H^{(k)})$ over $\T$, we have that
$$
\Nablae_V A = \pi \Nablae_V A.$$
Now we can simply apply formula (\ref{derivepi}) to obtain the
desired conclusion.
\eproof

\begin{proposition} \label{CDstarpr}
Suppose $f,g : \T \ra C^\infty(M)$ are smooth functions. Then
$T^{(k)}_fT^{(k)}_g$ is a smooth family of Toeplitz operators over
$\T$ and for any vector field $V$ on $\T$ we have that
\[\Nablae_V(T^{(k)}_fT^{(k)}_g) \sim \Nablae_V(T^{(k)}_{f\BTstar g}),\]
i.e. for all $L\in \bZ_+$ we have that
\[\left\| \Nablae_V (T^{(k)}_{f}T^{(k)}_g) - \sum_{l=0}^L \Nablae_V (T^{(k)}_{c_l(f,g)}
)\right\| = O(k^{-(L+1)}).\]
\end{proposition}

\proof

Using the notation of \cite{Sch} we let
\[A_L = D^L_\varphi T_f T_g - \sum_{l=0}^{L} D^{L-l}_\varphi T _{c_l(f,g)}.\]
Suppose $V$ be a vector field on $\T$. We will now establish by
induction that $\Nablae A_L$ is a zero order Toeplitz operator and
\[\s_0(\Nablae_V A_L ) = V[c_L(f,g)].\]
Note that Lemma \ref{sd=ds} implies this claim for
$A_0 = T_fT_g$, since we have just argued this operator is smooth.
Assume we have established this claim for $\Nablae A_{L-1}$.
Since
\[A_L = D_\varphi A_{L-1} - D_\varphi T_{c_{L-1}(f,g)},\]
we see that $A_L$ is a smooth family of Toeplitz operators
parameterized by $\T$ and
\[\Nablae_V A_L = D_\varphi \Nablae_V A_{L-1} - D_\varphi \Nablae_V T_ {c_{L-1}(f,g)}.\]
We see this is at most a first order operator by induction, but
\[\s_1(\Nablae_V A_L) = t ( \s_0(\Nablae_V A_{L-1}) - V[c_{L-1}(f,g)]) = 0\]
by the previous lemma, so it is at most a $0$-order operator. Applying the previous lemma
again we see that
\[\s_0(\Nablae_V A_L) = V[\s_0(A_L)] = V[c_L(f,g)].\]
This completes the inductive step.

The estimates of the theorem now follow since by induction $\Nablae_V
A_L$ is a zero order Toeplitz operator for all $L\in \bZ_+$.

\eproof

\begin{lemma}\label{Deriv}
The formal operator $D_V$ is a derivation for $\BTstar_\sigma$ for each $\s\in\T$, i.e.
\[D_V(f\BTstar g) = D_V(f)\BTstar g + f\BTstar D_V(g)\]
for all $f,g\in C^\infty(M)$.
\end{lemma}

\proof
By definition of $\BTstar_\s$ we have for all $\s\in \T$ that
\[T^{(k)}_{f,\s}T^{(k)}_{g,\s} \sim T^{(k)}_{f\BTstar_\s g,\s}.\]
Let $V$ be a vector field on $\T$.
By Proposition \ref{CDstarpr} we have that
\[\Nablae_V(T^{(k)}_{f}T^{(k)}_{g}) \sim \Nablae_V(T^{(k)}_{f\BTstar g}).\]
Considering the left hand side, we see that
\[
\Nablae_V(T^{(k)}_{f}T^{(k)}_{g})  =
   \Nablae_V(T^{(k)}_{f})T^{(k)}_{g} +  T^{(k)}_{f}\Nablae_V(T^{(k)}_{g}).
   \]
Now apply Theorem \ref{MainFGHCI} to get the wanted conclusion.
\eproof

\begin{proposition}\label{Curvature}
For two vector fields $V_1,V_2$ on $\T$, we have the formula
\begin{equation}
([\Nablae_{V_1},\Nablae_{V_2}] - \Nablae_{[V_1,V_2]}
)(T^{(k)}_f)  \sim T^{(k)}_{([D_{V_1},D_{V_2}]-
D_{[V_1,V_2]})(f)(1/(2k+n))}
\end{equation}
\end{proposition}

From this Proposition we conclude in particular that flatness of
$\Nablae$ implies flatness of $D$.

\proof
By \eqref{NablaeT} we see that
\begin{eqnarray*}
\Nablae_{V_1}(\Nablae_{V_2}(T_f^{(k)})) & = &
 \Nablae_{V_1}(T_{V_2[f]}^{(k)}) +
 \Nablae_{V_1}(T_{V'_2[F]f}^{(k)}) - \Nablae_{V_1}(T_{V'_2[F]}^{(k)}T_{f}^{(k)})\\
 & & + \frac1{2k+n}(\Nablae_{V_1}(T^{(k)}_{E(V_2)(f)}) -
 \Nablae_{V_1}(T_{H(V_2)}^{(k)}T_{f}^{(k)})).
\end{eqnarray*}
But now by applying first Proposition \ref{CDstarpr} and then
Theorem \ref{MainFGHCI} to this expression, we see that
\[\Nablae_{V_1}(\Nablae_{V_2}(T_f^{(k)})) \sim
T^{(k)}_{D_{V_1}D_{V_2}f}.\] The Proposition then follows from this
and Theorem \ref{MainFGHCI}. \eproof

\section{Formal trivializations and symmetry-invariant $*$-products}\label{PTSI}

\begin{definition}\label{formaltrivi2}
A formal trivialization of a formal connection $D$ is a
smooth map $P : \T \ra \D_h(M)$ which modulo $h$ is an isomorphism for all $\s\in \T$ and
such that
\[D_V(P(f)) = 0\]
for all vector fields $V$ on $\T$ and all $f\in C^\infty_h(M)$.
\end{definition}

Clearly if $D$ is not flat, such a formal trivialization will not exist even
locally on $\T$. However, if $D$ is flat, which is implied if
$\Nabla$ is projectively flat by Proposition \ref{Curvature}, then we have the
following result.

\begin{proposition}
Assume that $D$ is flat and that $\tD = 0$ mod $h$. Then locally around any point in $\T$
there exists a formal trivialization. If $H^1(\T,\bR) = 0$ then there exists a formal
trivialization defined globally on $\T$. If further $H^1_\Gamma(\T,\D(M)) = 0$
then we can construct $P$ such that it is $\Gamma$-equivariant.
\end{proposition}

In this proposition $H^1_\Gamma(\T,\D(M))$ simply refers to
the $\Gamma$-equivariant first de Rham cohomology of $\T$ with
coefficients in the real $\Gamma$-vector space $\D(M)$.

\proof
We write the formal trivialization we seek as
$$
P = \sum_{l=0}^\infty P_l h^l$$
where $P_l : \T \ra \D(M)$.
We need to solve
$$
D_V P = \sum_{l=0}^\infty V[P_l] h^l + \sum_{l=0}^\infty
\sum_{r=1}^{l} \tD^{(r)}(V)P_{l-r} h^l.$$
Hence we need that
\begin{equation}
V[P_l] = \sum_{r=1}^{l} \tD^{(r)}(V)P_{l-r}.\label{EqP}
\end{equation}
Now $P_0 = \Id$ solves this equation for $l=0$.
Assume that we have solved (locally, globally on $\T$ respectively
$\Gamma$-equivariantly on $\T$) this equation for $P_r$, for $r<l$. Then let
$\alpha_l \in \Omega^1(\T, \D(M))$ be given by
$$
\alpha_l(V) = \sum_{r=1}^{l} \tD^{(r)}(V)P_{l-r}.$$
We observe that $\alpha_l$ is $\Gamma$-invariant.
A short computation shows that the flatness of $D$ implies that $\alpha_l$
is closed on $\T$. Hence we can (locally, globally on $\T$ respectively
$\Gamma$-equivariantly on $\T$) solve (\ref{EqP}) for $P_l$.

\eproof

Now suppose we have a formal trivialization $P$ of the formal
 Hitchin connection $D$ determined by \eqref{Tdf}. Then $P$ is
 constant mod $h$ and we may and will assume that
\[P = \Id \mod h.\]
We can then define a new smooth family of star products,
parametrized by $\T$, by
\[f\star_\s g = P_\s^{-1}(P_\s(f) \BTstar_\s P_\s(g))\]
for all $f,g\in C^\infty(M)$ and all $\s\in \T$.

\begin{proposition}\label{ftrdq}
The star-products $\star_\s$ are independent of $\s\in\T$.
\end{proposition}

\proof Let $f,g \in C^\infty(M)$.
Since $D_V$ is a derivation for $\BTstar_\sigma$ for any $\s\in \T$, we have that
\[D_V(P_\s(f)\BTstar_\s P_\s(g)) = 0.\]
However
\begin{eqnarray*}
D_V(P_\s(f)\BTstar_\s P_\s(g)) & = & \tD(V)(P_\s(f)\BTstar_\s P_\s(g))\\
                         && + V[P_\s](f) \BTstar_\s P_\s(g)\\
                         && + P_\s(f) \BTstar_\s V[P_\s](g)\\
                         && + P_\s(f) V[\BTstar_\s] P_\s(g),
\end{eqnarray*}
which we compare with
\begin{eqnarray*}
V[f\star_\s g ]& = & V[P_\s^{-1}](P_\s(f) \BTstar_\s P_\s(g))\\
                         && + P_\s^{-1}(V[P_\s](f) \BTstar_\s P_\s(g))\\
                         && + P_\s^{-1}(P_\s(f) \BTstar_\s V[P_\s](g))\\
                         && + P_\s^{-1}(P_\s(f) V[\BTstar_\s] P_\s(g))
\end{eqnarray*}
and conclude, since $P V[P^{-1}] = \tD(V)$, that
\[V[f\star_\s g ] = 0.\]
\eproof

From the above we conclude Theorem \ref{general}.

Now, let us analyze equivalences between symmetry invariant $*$-products.
Suppose we have two differential $\Gamma$-invariant $*$-products $*$ and $*'$,
which are equivalent under some equivalence
$$ T = \Id + \sum_{j=1}^\infty h^j T_j,$$
where $T_j : C^\infty(M) \ra C^\infty(M)$ is a linear map for each $j\in \bN$, such that
$$T(f*g) = T(f)*'T(g).$$
By Theorem 2.22 in \cite{GR} it follows that $T_j$ is a differential operator for all $j\in\bN$.

\begin{proposition}
If the first discrete cohomology of $\Gamma$ with coefficient in the $\Gamma$-module
$C^\infty_0(M)$,
$H^1(\Gamma, C^\infty_0(M))$ and first de Rham cohomology of $M$ with real coefficients,
$H^1(M,\bR)$, both vanish, then we can find a $\Gamma$-invariant equivalence
between $*$ and $*'$.
\end{proposition}

\proof

We consider the given equivalence $T$ and by a short computation we get that
$$T^{-1} \gamma^* (T)(f*g) = T^{-1} \gamma^*(T)(f) * T^{-1} \gamma^*(T)(g),$$
hence $T^{-1}\gamma^*(T)$ is an automorphism of $*$. Since $H^1(M,\bR) = 0$
we get by proposition 3.3 in \cite{GR} that
there exists
$a_\gamma \in C^\infty_h(M)$ for each $\gamma\in \Gamma$ such that
\begin{equation}
T^{-1} \gamma^*(T) = \exp(\ad_* a_\gamma). \label{Tinner}
\end{equation}
We observe for $u,f\in C^\infty(M)$ that
\[\ad_*(u)(f) =  \{u,f\}h + O(h^2),\]
and that $\ad_*(c) =0 $ for all $c\in \bC_h$.
If we have that $a_\gamma = \sum_{j=0}^\infty a_\gamma^{(j)}h^j$,
then we may assume that $a^{(j)}_\gamma\in C^\infty_{0}(M)$. Furthermore $a_\gamma$ is
then uniquely determined by \eqref{Tinner}. Let us now assume that
\[a^{(i)}_\gamma = 0\]
for all $\gamma\in \Gamma$ and $i=0, \ldots j-1$. We
will then show that we can modify $T_{j}$ to obtain a
new equivalence which through \eqref{Tinner} produces a new
$a_\gamma$, which vanishes modulo $h^{j+1}$.

First we see that
\begin{eqnarray*}
\exp( \ad_* a_{\gamma_1\gamma_2}) & = & \exp (\ad_* a_{\gamma_1})
\exp(\ad_* \gamma_1^*(a_{\gamma_2}))\\
&= & \exp (\ad_* (a_{\gamma_1}\circ_* \gamma_1^*(a_{\gamma_2})))
\end{eqnarray*}
by Lemma 4.1 in \cite{GR}, where $\circ_*$ is the
Campbell-Baker-Hausdorff composition
\[a\circ_* b = a + \int_0^1 \psi(\exp (\ad_* a)\circ
\exp(t\ad_* b))b\quad dt,\]
where
\[\psi(z) = \frac{z \log(z)}{z-1}.\]
One has for $a,b\in C^\infty(M)$ that
\[a\circ_* b = a+b  + O(h^1).\]
From the above we may conclude that
\[a^{(j)}_{\gamma_1\gamma_2} = a^{(j)}_{\gamma_1} +  \gamma^*_1(a^{(j)}_{\gamma_2}).\]
Hence we see that
\[(a_\gamma^{(j)}) \in Z^1(\Gamma, C^\infty_0(M)).\]
But by assumption $H^1(\Gamma, C^\infty_0(M))= 0$, so this means
that there exists $a^{(j)}\in C^\infty_0(M)$ such that
\[a^{(j)}_\gamma = \gamma^*(a^{(j)}) - a^{(j)}.\]
Now replace $T_{j}$ by $T_{j} \exp(\ad_* (a^{(j)} h^j))$ and obtain a new equivalence
which produces a new $a_\gamma$ with the required vanishing. By induction we have the conclusion of the Proposition.
\eproof

\begin{remark}
From this we conclude that if the commutant of $\Gamma$ in $D(M)$ is
trivial, i.e. it contains only scalar multiples of the identity,
then a $\Gamma$ invariant differential $*$-product on $M$ is unique.
\end{remark}

\section{Applications to moduli space of flat connections}\label{sec5}

Let $\Sigma$ be a compact surface. Let $\cM$ be the moduli space
of flat $SU(n)$-connections on $\Sigma$
\[\cM = \Hom(\pi_1(\Sigma), SU(n))/SU(n).\]
There is a natural Poisson structure on $\cM$ (see \cite{FR1} and
\cite{FR2}). The symplectic leaves are specified by fixing the
conjugacy-class of the holonomy around each component of the
boundary of $\Sigma$. Let $(M,\omega)$ be a smooth symplectic leaf
of $\cM$. Pick a prequantum line bundle on $(M,\omega)$. Then the assumptions
of Theorem \ref{MainGHCI} are satisfied and we get the existence of the Hitchin
connection given by (\ref{Hitcon}).

As a corollary of Lemma \ref{Vriccipot}, we get that
\begin{theorem}
The Hitchin connection agrees with the connection constructed by
Axelrod, della Pietra and Witten in \cite{ADW}.
\end{theorem}

Now the techniques used in \cite{H} to show that the Hitchin connection
is projectively flat applies in our situation, hence we conclude that the
induced connection in the endomorphism bundle $\Nablae$ is flat. Since Teichm{\"u}ller
space is contractible, we can apply Theorem \ref{FQ} to obtain Theorem
\ref{MCGinvDG}.

\end{document}